\documentclass[12pt]{article}

\usepackage{amsthm}
\usepackage{amsmath}

\usepackage[usenames]{color}
\usepackage{amssymb}
\usepackage{amsmath}
\usepackage{amsthm}
\usepackage{amsfonts}
\usepackage{amscd}
\usepackage{graphicx}

\usepackage{fullpage}
\usepackage{float}
\usepackage{graphics}
\usepackage{latexsym}
\usepackage{epsf}
\usepackage[colorlinks=true,
urlcolor=red,
citecolor=webgreen,
linkcolor=blue]{hyperref}
\definecolor{webgreen}{rgb}{0.0,0.5,0.0}

\newcommand{\seqnum}[1]{\href{https://oeis.org/#1}{\rm \underline{#1}}}

\DeclareMathOperator{\mathSS}{SS}
\DeclareMathOperator{\SCS}{SCS}
\DeclareMathOperator{\SSCS}{SSCS}

\theoremstyle{plain}
\newtheorem{theorem}{Theorem}

\newtheorem{proposition}[theorem]{Proposition}

\theoremstyle{definition}
\newtheorem{definition}[theorem]{Definition}

\theoremstyle{remark}

\begin{document}

\title{Counting Various Classes of \\  Tournament Score Sequences}       
\author{Paul K. Stockmeyer \\ \ \\
         College of William \& Mary \\
         Williamsburg, Virginia, USA \\
       Email: \href{mailto:stockmeyer@cs.wm.edu}{stockmeyer@cs.wm.edu}}       
\date{}        
\maketitle

\begin{abstract}  Based on known methods for computing the number of distinct score sequences for $n$-vertex tournaments,
we develop algorithms for computing the number of distinct score sequences for self-com\-ple\-men\-ta\-ry tournaments, strong
tournaments, and tournaments that are both self-complementary and strong.
\end{abstract}

\section{Introduction}\label{intro}

In graph theory, a \emph{tournament}, sometimes called a \emph{round-robin 
tournament}, is an oriented complete graph, i.e., a directed graph where for 
any two vertices $v$ and $w$ there is either a directed edge from $v$ to $w$ 
or from $w$ to $v$ but not both. For background information about tournaments,
see the classic book by Moon \cite{moon} or the article by Harary and Moser 
\cite{HandM}.

If an edge of a tournament is directed from vertex $v$ to vertex $w$ we say that 
$v$ \emph{dominates} $w$.  The \emph{score} of a vertex $v$ is the number of 
vertices that $v$ dominates.  The \emph{score sequence} of an $n$-vertex tournament 
is the sequence $(s_1, s_2, \ldots , s_n)$ of vertex scores, written in 
nondecreasing order:
\begin{equation}\label{nondecreasing}
     s_1 \le s_2 \le \cdots \le s_n.
\end{equation}
For example, we observe that the score sequences for 5-vertex 
tournaments are $(0,1,2,3,4)$, $(0,1,3,3,3)$, $(0,2,2,2,4)$, $(0,2,2,3,3)$, 
$(1,1,1,3,4)$, $(1,1,2,2,4)$, $(1,1,2,3,3)$, $(1,2,2,2,3)$, and $(2,2,2,2,2)$.

This paper develops methods for computing the number of score sequences for 
various classes of $n$-vertex tournaments.  Section \ref{history} reviews the 
history of computing the total number of $n$-vertex score sequences, introducing tools 
that will be used in later sections. Section \ref{selfcomp} provides an algorithm 
for counting self-complementary score sequences, Section \ref{strong} contains 
two different methods for counting strong score sequences, while Section 
\ref{strongselfcomp} gives two methods for counting score sequences that are 
both strong and self-complementary.  Section \ref{observations} contains some 
asymptotic results and conjectures.

\section{Counting score sequences}\label{history}

Sequence \seqnum{A000571} in the {\it On-Line Encyclopedia of Integer Sequences} 
(OEIS) \cite{oeis} gives the number of different score sequences that are possible 
for an $n$-vertex tournament. Postulating the existence of a 0-vertex tournament 
with null score sequence, this sequence begins with $n=0$ and reads 1, 1, 1, 2, 4, 
9, 22, 59, 167, 490, 1486, \dots\ .  There is no known simple formula for this 
sequence but various methods have been devised for determining its terms.

The first major effort to catalog score sequences was carried out in 1920 by 
MacMahon \cite{mac}, who expanded the product
$$
\prod_{1\le i<j \le n}\left(a_i + a_j\right)
$$
to generate the score sequence of each of the $2^{\binom{n}{2}}$ labeled tournaments 
on $n$ vertices.  For example, with $n=5$, fourteen of the 1024 terms in this product 
are $a_1^2a_2^2a_3^1a_4^3a_5^2$, representing some of the labeled tournaments with 
score sequence $(1,2,2,2,3)$.  By combining these terms with the other terms with 
the same collection of exponents, MacMahon calculated that there were $14\times20= 280$ 
labeled tournaments with this score sequence, the most of any of the nine possible 
score sequences.  In a Herculean effort of hand calculation, MacMahon carried out 
this procedure up through $n=9$, determining the number of labeled tournaments 
possessed by each of 490 possible score sequences that arose.  The most common 
sequence turned out to be $(2,3,3,4,4,4,5,5,6)$, occurring 5,329,376,640 times 
among the $2^{36} =  68,719,476,736$ labeled tournaments on nine vertices.

Presumably unaware of the work of MacMahon, David \cite{david} made essentially 
the same hand calculations in 1959, but only for $n\le 8$.  This was followed by 
Alway \cite{alway}, who in 1962 used a computer to catalog all score sequences 
up to $n = 10$.

Of course it is not necessary to generate all the $n$-vertex tournament score 
sequences in order to count them.  The following test for a sequence of non-negative 
integers to be the score sequence of a tournament is usually attributed to 
Landau \cite{landau}. 

\begin{proposition}
A sequence $(s_1, s_2, \dots , s_n)$ of $n \ge 1$ non-negative integers satisfying \eqref{nondecreasing}
is the score sequence of a tournament if and only if
\begin{align}
   &\sum_{i=1}^r s_i \ge \binom{r}{2}\mbox{ for }1\le r < n, \mbox{ and}\label{rule1}\\
   &\sum_{i=1}^n s_i = \binom{n}{2}.\label{rule2}
\end{align}
\end{proposition}
This test was utilized by Bent \cite{bent} in a remarkable 1964 Master of Science 
dissertation, also published by Narayana and Bent \cite{NandBent}.  We present 
Bent's algorithm, which will be used in later sections.

\begin{definition}
For all positive $n$ and all non-negative $T$ and $E$, the array entry 
$F_n\left[T,E\right]$ is the number of sequences of length $n$ satisfying 
(\ref{nondecreasing}), (\ref{rule1}), $s_n=E$, and
$$
                         \sum_{i=1}^n s_i = T.
$$
\end{definition}

\noindent
(We can imagine that $T$ stands for {\it Total} and $E$ for {\it End}.) These 
$F_n$ arrays can be computed using the following recursive formulas.

\begin{proposition}[Bent]\label{bent0}
We have
\begin{equation*}
     F_1\left[T,E\right] = \begin{cases}1,&\text{if $T=E$};\\0,&\text{otherwise},\end{cases}
\end{equation*}
and for $n\ge 2$ we have
\begin{equation*}
        F_n\left[T,E\right] = \begin{cases}
                  \displaystyle{\sum_{k=0}^E F_{n-1}\left[T\!-\!E,k\right]},
                  &\text{if $T\!-\!E\ge \binom{n-1}{2}$};\\
                       0,&\text{otherwise}.\end{cases}
\end{equation*}
\end{proposition}

\noindent
Bent then used the values in the $F_n$  arrays to count score sequences.

\begin{theorem}[Bent] \label{bent3} Let $S(n)$ denote the number of distinct score 
sequences for $n$-vertex tournaments.  Then for $n\ge 1$ we have
\begin{equation*}
   S(n) = \sum_{E=\left\lceil\frac{n-1}{2}\right\rceil}^{n-1}F_n\left[\binom{n}{2},E\right].
\end{equation*}
\end{theorem}

Using Theorem \ref{bent3} and the computer available to him, Bent computed $S(n)$ 
for $n\le 27$.  With an alternative method, which allowed him to compute $S_{2n}$ 
from array $F_n$  and to compute $S_{2m+1}$ from arrays $F_m$ and $F_{m+1}$, he was 
able to extend his result up to $n=36$.  The limiting factor was the speed of the 
computer he was using. 

A recent paper by Claesson et al.\ \cite{CDFS} provides
an alternative method for counting score sequences.
Today the OEIS entry \seqnum{A000571} gives values of $S(n)$ 
up to $n = 1,675$.  As a test of our programs for this paper, we verified these 
numbers for $n \le 500$.

\section{Counting self-complementary score sequences}\label{selfcomp}

The {\it complement} of a tournament $T_n$ is the tournament $T_n^c$ obtained from 
$T_n$ by reversing the direction of all of its edges. Thus a vertex with score $s_i$ 
in $T_n$ becomes a vertex of score $n\!-\!1\!-\!s_i$ in $T_n^c$.  A tournament is 
called {\it self-complementary} if it is isomorphic to its complement.  In this 
case, for each vertex of score $s_i$ there is a corresponding vertex of score
$n\!-\!1\!-\!s_i$.

\begin{definition}\label{self}
A score sequence of length $n$ is called {\it self-complementary} if $s_{n+1-i} 
= n\!-\!1\!-\!s_i$ for $1\le i \le \left\lfloor\frac{n}{2}\right\rfloor$.
\end{definition}

\noindent
For example, the self-complementary score sequences of length $n=6$ are 
$(0,1,2,3,4,5)$, $(0,2,2,3,3,5)$, $(1,1,1,4,4,4)$, $(1,1,2,3,4,4)$, $(1,2,2,3,3,4)$, 
and $(2,2,2,3,3,3)$.  Clearly self-complementary tournaments have self-complementary 
score sequences.  The converse, however, is not true.  Of the four non-isomorphic 
tournaments with score sequence $(1,1,2,3,4,4)$, two are self-complementary and 
the other two are complements of each other.

Sequence \seqnum{A345470} in the OEIS, authored by Givner in 2021, lists 
the number of self-complementary score sequences for $n$-vertex tournaments.  
Givner wrote a program that generated all length $n$ score sequences for $n\le 34$ 
and then counted those that were self-complementary.  As with Bent's program, 
computer time was the limiting factor in Givner's calculations.

Without actually generating the self-complementary score sequences, we can count 
them using the $F_n$ arrays of Bent from Section \ref{history}. 
 
\begin{theorem}\label{selfcom}
Let $\SCS(n)$ denote the number of distinct self-complementary score sequences of 
length $n$.  Then for all $m\ge 1$ we have
\begin{equation*}
   \SCS(2m) = \sum_{T=\binom{m}{2}}^{m(m-1)}\ \sum_{E=\left\lceil \frac{T}{m}
   \right\rceil}^{m-1}F_m\left[T, E\right],
\end{equation*}
and
\begin{equation*}
   \SCS(2m+1) = \sum_{T=\binom{m}{2}}^{m^2}\ \sum_{E=\left\lceil \frac{T}{m}
   \right\rceil}^{m}F_m\left[T, E\right].
\end{equation*}
\end{theorem}

\begin{proof}
From Definition \ref{self} we know that a self-complementary score sequence of 
length $2m$ is determined by its first $m$ terms. Therefore $\SCS(2m)$ is a certain 
sum of values in the array $F_m$.  For the lower summation limit on $E$ we note 
that $E=s_m$ must be at least as large as the average of $s_1$ through $s_m$, or 
$\frac{T}{m}$.  For the upper limit on $E$ we observe that $s_m \le s_{m+1}$ and 
$s_m + s_{m+1} = 2m-1$.  This implies that $E=s_m \le \left\lfloor\frac{2m-1}{2}
\right\rfloor = m-1$.  The lower summation limit on $T$ comes from equation 
(\ref{rule1}) while for the upper limit we note that $T\le m\cdot s_m \le m(m-1)$.

For self-complementary score sequences of length $2m+1$ we have that $s_{m+1}$ 
must be $m$, and the score sequence is again determined by its first $m$ terms.  
The only difference between this case and the even case is that $E=s_m$ can now 
be as large as $s_{m+1}=m$, and $T$ can now be as large as $m\cdot s_m \le m^2$.
\end{proof}

We used Theorem \ref{selfcom} to compute $\SCS(n)$ for $n$ from 2 to 500 as 
displayed in the OEIS entry \seqnum{A345470}, confirming the results of Givner 
for $n$ from 2 to 34. 

\section{Counting strong score sequences}\label{strong}

A tournament is called {\it strong,} or {\it strongly connected,} if every vertex 
of the tournament can reach every other vertex along a directed path.  A tournament 
is {\it reducible} if its vertex set can be partitioned into two nonempty sets
$A$ and $B$ with every vertex in $B$ dominating every vertex in $A$.  It is well 
known that a tournament is strong if and only if it is not reducible.

Harary and Moser \cite[Theorem 9]{HandM} showed that the property of being strong 
can be determined from a tournament's score sequence.

\begin{proposition}\label{strongseq}
A tournament is strong if and only if its score sequence satisfies equations 
$(\ref{nondecreasing})$, $(\ref{rule2})$ and
\begin{equation}\label{strongrule}
   \sum_{i=1}^r s_i > \binom{r}{2}\text{ for }1 \le r < n.
\end{equation}
\end{proposition}

We call score sequences that satisfy equation (\ref{strongrule}) {\it strong} score 
sequences. For example, the seven strong score sequences for $n=6$ are 
$(1,1,2,3,4,4)$, $(1,1,3,3,3,4)$, $(1,2,2,2,4,4)$, $(1,2,2,3,3,4)$, $(1,2,3,3,3,3)$, 
$(2,2,2,2,3,4)$, and $(2,2,2,3,3,3)$.  

In this section we present two algorithms for counting strong score sequences.  Our 
first method is a modified version of Bent's algorithm of Section \ref{history}, 
using $G_n$ arrays based on equation (\ref{strongrule}) instead of the $F_n$ arrays 
based on equation (\ref{rule1}).

\begin{definition}\label{defG}
For all positive $n$ and all non-negative $T$ and $E$, let $G_n\left[T,E\right]$  
be the number of sequences of length $n$ satisfying (\ref{nondecreasing}), 
(\ref{strongrule}), $s_n=E$, and
$$
                         \sum_{i=1}^n s_i =T.
$$
\end{definition}

\begin{proposition}\label{G0}
We have
\begin{equation*}
     G_1\left[T,E\right] = \begin{cases}1,& \text{if $T=E$};\\0,&\text{otherwise},\end{cases}
\end{equation*}
and for $n\ge 2$ we have
\begin{equation*}
        G_n\left[T,E\right] = \begin{cases}
                  \displaystyle{\sum_{k=\left\lceil\frac{T-E}{n-1}\right\rceil}^E 
                   G_{n-1}\left[T\!-\!E,k\right]},
              &\text{if $T\!-\!E > \binom{n-1}{2}$};\\
                       \quad 0,&\text{otherwise}.\end{cases}
\end{equation*}
\end{proposition}

\begin{proof}
The case $n=1$ is clear: a sequence of length one gets counted if and only if its 
last (and only) element $E$ is equal to the sum $T$ of all its elements.  For 
$n\ge 2$, consider a sequence $S = (s_1, \ldots ,s_n)$ counted by $G_n(T,E)$.
We have $\sum_{i=1}^r s_i > \binom{r}{2}$ for $r<n$, $\sum_{i=1}^n s_i = T$, and 
$s_n = E$.  If we delete the last term $s_n$ we have a sequence $\widehat{S} = 
(s_1,\ldots, s_{n-1})$ with $\sum_{i=1}^{n-1} s_i = T\!-\!E$ and $s_{n-1}$ equal 
to some integer $k$ with $\frac{T-E}{n-1} \le k \le E$. The limits arise because 
$s_{n-1}$ must be at least as large as  of the average of the first $n-1$ terms 
of the sequence and no larger than $s_n$.  Thus for every sequence $S$ counted 
by array element $G_n\left[T,E\right]$, the corresponding sequence $\widehat{S}$ 
is counted  by one of the indicated $G_{n-1}\left[T\!-\!E,k\right]$ array elements.  
However, not all strings counted by these array elements satisfy the requirement 
that $\sum_{i=1}^{n-1} s_i > \binom{n-1}{2}$; our sum only includes those terms 
for which this requirement is satisfied.
\end{proof}

Our first method for counting strong score sequences now follows. 

\begin{theorem}\label{strontheorem1}
Let $\mathSS(n)$ denote the number of distinct strong score sequences for $n$-vertex 
tournaments.  Then for  $n\ge 1$ we have
$$
  \mathSS(n) = \sum_{E=\left\lceil\frac{n-1}{2}\right\rceil}^{n-2}G_n\left[\binom{n}{2},E\right].
$$
\end{theorem}

\begin{proof}
By definition, $G_n\left[\binom{n}{2},E\right]$ counts strong score sequences of 
length $n$ with largest score being $E$.  This score is bounded below by the average 
score, and bounded above by $n-2$ for strong tournaments.
\end{proof}

Our second method for counting strong score sequences is to count reducible score 
sequences and subtract that number from the number of all score sequences.  Every 
reducible tournament can be characterized by a non-empty vertex set $A$ which induces 
a strong subtournament, and a nonempty set $B$ of the remaining vertices, with every 
vertex in $B$ dominating every vertex in $A$.  The score sequence of such a tournament 
consists of an initial strong score sequence of length $\left|A\right|$, followed 
by an arbitrary valid score sequence of length $\left|B\right|$ with each score 
increased by $\left|A\right|$.  Some typical reducible score sequences of length 
6 are $(0,2,3,3,3,4)$ with $\left|A\right| = 1$, $(1,1,1,3,4,5)$ with $\left|A\right|=3$, 
$(1,1,2,2,4,5)$ with $\left|A\right|=4$, and $(1,2,2,2,3,5)$ with $\left|A\right|= 5$.
This characterization of reducible score sequences provides us with the following 
bootstrapping method for counting strong score sequences, using only the series 
$S(n)$ that counts all score sequences. 

\begin{theorem}\label{strong2}
Terms in the sequence $\mathSS(n)$ with $n\ge 1$ can be computed recursively from the 
sequence $S(n)$ using the formula
$$
   \mathSS(n) = S(n) - \sum_{i=1}^{n-1}\mathSS(i)S(n-i).
$$
\end{theorem}

Claesson et al.\ \cite{CDFS} observed that this method also works when using their method for
computing $S(n)$.

We computed $\mathSS(n)$ for $n$ from 1 to 500 by both of these methods with identical 
results, as displayed in the OEIS entry \seqnum{A351822}. 

\section{Counting strong self-complementary score sequences}\label{strongselfcomp}

In this section we sketch two methods for computing the number of score sequences 
for $n$-vertex tournaments that are both strong and self-complementary.  The 
first is based on our method of counting self-complementary score sequences, but 
using the $G_n$ arrays rather than the $F_n$ arrays.

\begin{theorem}
Let $\SSCS(n)$ denote the number of distinct score sequences of length $n$ that are 
both strong and self-complementary. Then for all $m\ge 1$ we have
\begin{equation*}
   \SSCS(2m) = \sum_{T=\binom{m}{2}+1}^{m(m-1)}\ \sum_{E=\left\lceil \frac{T}{m}
                 \right\rceil}^{m-1}G_m\left[T, E\right]
\end{equation*}
and
\begin{equation*}
   \SSCS(2m+1) = \sum_{T=\binom{m}{2}+1}^{m^2}\ \sum_{E=\left\lceil \frac{T}{m}
             \right\rceil}^{m}G_m\left[T, E\right].
\end{equation*}
\end{theorem}

Our second method for counting strong self-complementary score sequences is similar 
to our second method for counting strong score sequences:  subtract the number 
of self-complementary reducible score sequences from the number of self-complementary
score sequences.  Now a self-complementary reducible score sequence of length $n$ 
begins with a strong score sequence of length $i$ with $1 \le i \le \lfloor\frac{n}{2}
\rfloor$ and ends with the complement of this initial subsequence. In between we 
must have a self-complementary sequence of length $n-2i$, with each term increased 
by $i$.  For example, the three self-complementary reducible score sequences of 
length $n=6$ are $(0, 1, 2, 3, 4, 5)$ and $(0, 2, 2, 3, 3, 5)$ with $i=1$, and 
$(1, 1, 1, 4, 4, 4)$ with $i=3$.  Then the sequence that counts strong self-complementary 
score sequences can be computed from the sequences for strong sequences and 
self-complementary sequences.

\begin{theorem}  For all $n\ge 1$ we have 
$$
  \SSCS(n) =  \SCS(n) - \sum_{i=1}^{\lfloor \frac{n}{2}\rfloor} \mathSS(i)\SCS(n-2i).
$$
\end{theorem}

We observe that for $n\!=\!2m$ and $i\!=\!m$ the last term in the above summation 
is $\mathSS(m)\SCS(0)$.  For this to be correct we must have $\SCS(0) = 1$, i.e., we must 
assume the existence of a self-complementary tournament on zero vertices. 

We computed $\SSCS(n)$ for $n$ from 1 to 500 by both of these methods with identical 
results, as displayed in the OEIS entry \seqnum{A351869}. 

\section{Asymptotic Results and Conjectures}\label{observations}

Moon \cite{moon} reports that Erd\H{o}s and Moser, in an unpublished work, 
showed that there exist constants $c_1$ and $c_2$ such that
\begin{equation*}
    \frac{c_1\cdot 4^n}{n^{9/2}} < S(n) < \frac{c_2\cdot 4^n}{n^{3/2}}.
\end{equation*}
Winston and Kleitman \cite{winston} improved these bounds, showing that
$$
    \frac{c_1\cdot 4^n}{n^{5/2}} < S(n) < \frac{c_2\cdot 4^n}{n^{2}}
$$
and conjecturing that
$$
   S(n) =\Theta\left(\frac{4^n}{n^{5/2}}\right).
$$
This conjecture was confirmed by Kim and Pittel \cite{kim}.  Vaclav Kotesovec 
asserts in the OEIS entry \seqnum{A000571} that
$$
  S(n) \sim \frac{c \cdot 4^n }{ n^{5/2}},
$$
with $c = 0.392478 \dots\ $, but with no apparent proof.

For strong score sequences, Kotesovec asserts in \seqnum{A351822} that
$$
  \mathSS(n) \sim \frac{c \cdot 4^n }{ n^{5/2}},
$$
with $c = 0.202756 \dots\ $.  It would be nice to at least have a proof of the 
order of magnitude of $\mathSS(n)$.  Assuming that the assertions of Kotesovec are 
correct, we have that slightly more than half of all tournament score sequences 
are strong.  

As for self-complementary score sequences, we conjecture that both
$$
 \SCS(n) = \Theta\left(\frac{2^n}{n^{3/4}}\right) \text{ and } 
        \SSCS(n) = \Theta\left(\frac{2^n}{n^{3/4}}\right).
$$
We have no theoretical evidence for these formulas, but they are consistent with the data available.
Results suggest that over seventy percent of all self-complementary score 
sequences are strong.

\bigskip\hrule\bigskip

\noindent
2020 \emph{Mathematics Subject Classification}: Primary 05C30; Secondary 05A17,	05C20.

\medskip

\noindent
\emph{Keywords}: tournament, score sequence, strong, self-complementary.

\bigskip\hrule\bigskip

\noindent (Concerned with OEIS sequences \seqnum{A000571}, \seqnum{A345470}, \seqnum{A351822}, and \seqnum{A351869}.)

\end{document}